\documentclass[11pt]{article}

\usepackage{epsfig}
\usepackage{amssymb, latexsym}
\usepackage{amscd}
\usepackage[all]{xy}
\usepackage{epsf}
\usepackage[mathscr]{eucal}

\usepackage{amsmath,amsfonts,amscd,amssymb,amsthm}

\usepackage{appendix}

\long\def\comment#1\endcomment{}

\theoremstyle{plain}

\newtheorem{theorem}{{\sc Theorem}}[section]
\newtheorem{lemma}[theorem]{\sc Lemma}

\newtheorem{prop}[theorem]{\sc Proposition}
\newtheorem{coroll}[theorem]{\sc Corollary}
\newtheorem{klemma}[theorem]{\sc Key-lemma}

\theoremstyle{plain}

\newtheorem{defn}[theorem]{\sc Definition}

\theoremstyle{exercise}
\newtheorem{remark}[theorem]{\sc Remark}

\newtheorem{example}[theorem]{\sc Example}

\makeatletter \@addtoreset{equation}{section} \makeatother

\def\eqref#1{\thetag{\ref{#1}}}

\let\latexref=\ref
\def\ref#1{{\normalfont{\latexref{#1}}}}

\newcommand{\mb}{{\:\raisebox{3pt}{\text{\circle*{1.5}}}}}

\def\dlim_#1{{\displaystyle\lim_{#1}}^\hdot}

\newcommand{\Ext}{\operatorname{Ext}}

\newcommand{\Ob}{\mathrm{Ob}}

\newcommand{\Tetra}{\mathscr{T}etra}

\newcommand{\opp}{\mathrm{opp}}
\newcommand{\BBar}{\mathrm{Bar}}
\newcommand{\Cobar}{\mathrm{Cobar}}
\newcommand{\GS}{\mathrm{GS}}
\newcommand{\Hom}{\mathrm{Hom}}
\newcommand{\Bimod}{\mathscr{B}imod}
\newcommand{\Bicomod}{\mathscr{B}icomod}
\newcommand{\Ind}{\mathrm{Ind}}
\newcommand{\Coind}{\mathrm{Coind}}

\newcommand{\RHom}{\mathrm{RHom}}

\def\wtilde#1{\widetilde{#1}\vphantom{#1}}

\title{\sc{Tetramodules over a bialgebra form a 2-fold monoidal category}}

\author{\sc{Boris Shoikhet}}
\date{}

\begin{document}\maketitle
{\footnotesize
\begin{center}{\parbox{4,5in}{{\sc Abstract.} Let $B$ be an associative bialgebra over any field.
A {\it module} over $B$ in the sense of deformation theory is a {\it tetramodule} over $B$. All tetramodules form an abelian category.
This category was studied by R.Taillefer [Tai1,2]. In particular, she proved that for any bialgebra $B$, the abelian category $\Tetra(B)$ has enough injectives, and that $\Ext^\mb(B,B)$ in this category coincides with the Gerstenhaber-Schack cohomology of $B$.

We prove that the category $\Tetra(B)$ of tetramodules over any bialgebra $B$ is a 2-fold-monoidal category, with $B$ a unit object in it.
Roughly, this means that the category $\Tetra(B)$ admits two monoidal structures, with common unit $B$, which are compatible in some rather non-trivial way (the concept of an $n$-fold monoidal category is introduced in [BFSV]).
Within (yet unproven) 2-fold monoidal analogue of the Deligne conjecture, our result would imply that $\RHom^\mb(B,B)$ in the category of tetramodules is naturally a {\it homotopy 3-algebra.}
 }}
\end{center}
}

\bigskip
\bigskip

\section*{\sc Introduction}
This paper contains a part of results of our 2009 archive preprint [Sh1], presented at a glance suitable for a journal publication. It covers a construction of a 2-fold monoidal structure on the category of tetramodules, with  all necessary definitions, and an overview of the results of R.Taillefer [Tai1,2] on tetramodules and the Gerstenhaber-Schack cohomology [GS] (formerly served as Appendix in [Sh1]), as well as a computation of the Gerstenhaber-Schack cohomology for the
free commutative cocommutative bialgebra $S(V)$, for a $V$ is a vector space.

Our approach to the $n$-fold monoidal Deligne conjecture from [Sh1] still remains unfinished, and is not considered here.

The paper is organized as follows.

In Section 1 we recall, in some detail, the definition of $n$-fold monoidal categories from [BFSV].

In Section 2 we define the abelian category of tetramodules over a bialgebra, and construct a 2-fold monoidal category structure on it.
This Section is the core of the paper.

In Section 3 we recall (with proofs) some of results of R.Taillefer [Tai1,2], which give a relation between our results in Section 2, and deformation theory of associative bialgebras. We slightly modify the original proof in presentation (with the same main ideas), which, we hope, makes it more readable. As an application, we compute the Gerstenhaber-Schack cohomology of the free commutative cocommutative bialgebra $S(V)$ (the answer is well-known to specialists, but we were unable to find any published proof of it).

The main application of the results of Section 2 claims, modulo the Deligne conjecture for $n$-fold monoidal abelian categories, that
$\RHom^\mb_{\Tetra(A)}(A,A)$ is a homotopy 3-algebra, when the bialgebra $A$ is a Hopf algebra (over any field).
However, the $n$-fold monoidal Deligne conjecture is unproven yet (for $n>1$). In [Sh1], we tried to use some constructions of [Sch] to prove it, but (so far) we did not succeed in this direction.

We suggested a new approach to the $n$-fold monoidal Deligne conjecture, based on the Kock-To\"{e}n's simplicial Deligne conjecture [KT], and on the Leinster's definition [L] of weak Segal monoids.
Our recent paper [Sh2] is a first step in realization of this approach. We hope to complete it in our next papers.

\subsection*{}
\subsubsection*{\sc Acknowledgements}
I am thankful to Tom Leinster and Dima Tamarkin for discussions on $n$-fold monoidal categories, and to Bernhard Keller for explanation of the formalism of $(P,Q)$-pairs (used in Section 3). Dima Kaledin explained to me what is the Deligne conjecture for monoidal abelian categories, which inspired many of constructions presented here. Alexey Davydov communicated to me on the S.Schwede's paper [Sch] several years before I started this work; my studying of this paper, as well as many discussions with Stefan Schwede himself, had being a great source of inspiration.

The work was done at the Universit\'{e} du Luxembourg and at the Max-Planck Insitut f\"{u}r Mathematik, Bonn.
I am thankful for these two Institutions for excellent working conditions and for financial support.
The work in Luxembourg was supported by research grant R1F105L15 of Martin Schlichenmaier.

\section{\sc $n$-fold monoidal categories}
\subsection{\sc Introduction}
An $n$-fold monoidal category $\mathcal{C}$ is a category with $n$
monoidal structures $\otimes_1,\dots,\otimes_n\colon
\mathcal{C}\times\mathcal{C}\to\mathcal{C}$ which obey some
compatibility relations, which may seem complicated when are written down explicitly. They were introduced in [BFSV], with the
following motivation.

When $\mathscr{C}$ is a (1-)monoidal category, its nerve $N\mathcal{C}$ a simplicial $H$-space, up to a homotopy.
In fact, all homotopies between (the higher) homotopies can be chosen in a coherent way, such that the classifying space $B\mathcal{C}$ becomes an algebra over the topological Stasheff operad. It implies that the group completion of $B\mathcal{C}$ is a 1-fold loop space, $B\mathcal{C}=\Omega X$ for some connected space $X$.

The $n$-fold loop spaces $\Omega^n(X)$ are defined iteratively, $\Omega^n(X)=\Omega(\Omega^{n-1}(X))$. The $n$-fold monoidal categories are defined iteratively, in a similar way. The definition achieves the following goal: the group completion of the classifying space $B\mathcal{C}$, for an $n$-fold monoidal category $\mathcal{C}$, is an $n$-fold loop space.

\subsection{\sc Definition}
\begin{defn}{\rm
A (strict) monoidal category is a category $\mathcal{C}$ together
with a functor $\otimes\colon
\mathcal{C}\times\mathcal{C}\to\mathcal{C}$ and an object
$\mathbb{I}\in\Ob(\mathcal{C})$ such that
\begin{itemize}
\item[1.] $\otimes$ is strictly associative;
\item[2.] $\mathbb{I}$ is a strict two-sided unit for $\otimes$.
\end{itemize}
A monoidal functor $(F,\eta)\colon\mathcal{C}\to\mathcal{D}$ between
monoidal categories is a functor $F$ such that
$F(\mathbb{I}_\mathcal{C})=\mathbb{I}_\mathcal{D}$ with a natural transformation
\begin{equation}\label{eta}
\eta_{A,B}\colon F(A)\otimes F(B)\to F(A\otimes B)
\end{equation}
which satisfies the following conditions:
\begin{itemize}
\item[1.] Internal associativity: the following diagram commutes \begin{equation}\label{intassoc}
\xymatrix{ F(A)\otimes F(B)\otimes F(C)\ar[rr]^{\eta_{A,B}\otimes
id_{F(C)}}\ar[d]^{id_{F(A)}\otimes\eta_{B,C}}&&F(A\otimes B)\otimes
F(C)\ar[d]^{\eta_{A\otimes B,C}}\\
F(A)\otimes F(B\otimes C)\ar[rr]^{\eta_{A,B\otimes C}}&&F(A\otimes
B\otimes C) }
\end{equation}
\item[2.] Internal unit conditions:
$\eta_{A,\mathbb{I}}=\eta_{\mathbb{I},A}=id_{F(A)}$.
\end{itemize}
}
\end{defn}
The crucial in this definition is that {\it the map $\eta$ is not
required to be an isomorphism.}

Denote by ${\mathscr{MC}at}$ the category of (small) monoidal
categories and monoidal functors.
\subsubsection{\sc }
\begin{defn}
{\rm
A 2-fold monoidal category is a monoid in ${\mathscr{MC}at}$. This
means, that we are given a monoidal category
$(\mathcal{C},\otimes_1,\mathbb{I})$, and a monoidal functor
$(\otimes_2,\eta)\colon\mathcal{C}\times\mathcal{C}\to\mathcal{C}$
which satisfies the following axioms:
\begin{itemize}
\item[1.] External associativity: the following diagram commutes in
${\mathscr{MC}at}$
\begin{equation}\label{external}
\xymatrix{
\mathcal{C}\times\mathcal{C}\times\mathcal{C}\ar[rr]^{(\otimes_2,\eta)\times
id_{\mathcal{C}}}\ar[d]^{id_{\mathcal{C}}\times(\otimes_2,\eta)}&&\mathcal{C}\times\mathcal{C}\ar[d]^{(\otimes_2,\eta)}\\
\mathcal{C}\times\mathcal{C}\ar[rr]^{(\otimes_2,\eta)}&&\mathcal{C}}
\end{equation}
\item[2.] External unit conditions:
the following diagram commutes in ${\mathscr{MC}at}$
\begin{equation}\label{externalunit}
\xymatrix{
\mathcal{C}\times\mathbb{I}\ar[r]^{\subseteq}\ar[d]^{\cong}&\mathcal{C}\times\mathcal{C}\ar[d]^{(\otimes_2,\eta)}&\mathbb{I}\times\mathcal{C}
\ar[l]_{\supseteq}\ar[d]^{\cong}\\
\mathcal{C}\ar[r]^{=}&\mathcal{C}&\mathcal{C}\ar[l]_{=}}
\end{equation}
\end{itemize}
}
\end{defn}
Note that the role of the monoidal structures $\otimes_1$ and
$\otimes_2$ in this definition is not symmetric.

Explicitly the definition above means that we have an operation
$\otimes_2$ with the two-sided unit $\mathbb{I}$ (the same that for
$\otimes_1$) and a natural transformation
\begin{equation}\label{natural}
\eta_{A,B,C,D}\colon (A\otimes_2 B)\otimes_1(C\otimes_2 D)\to
(A\otimes_1 C)\otimes_2 (B\otimes_1 D)
\end{equation}
The internal unit conditions are:
$\eta_{A,B,\mathbb{I},\mathbb{I}}=\eta_{\mathbb{I},\mathbb{I},A,B}=id_{A\otimes_2B}$, and
the external unit conditions are:
$\eta_{A,\mathbb{I},B,\mathbb{I}}=\eta_{\mathbb{I},A,\mathbb{I},B}=id_{A\otimes_1B}$. As
well, one has the morphisms
\begin{equation}\label{12left}
\eta_{A,\mathbb{I},\mathbb{I},B}\colon A\otimes_1B\to A\otimes_2B
\end{equation}
and
\begin{equation}\label{12right}
\eta_{\mathbb{I},A,B,\mathbb{I}}\colon A\otimes_1B\to B\otimes_2A
\end{equation}

The internal associativity gives the commutative diagram:
\begin{equation}\label{intexpl}
\xymatrix{
(U\otimes_2V)\otimes_1(W\otimes_2X)\otimes_1(Y\otimes_2Z)\ar[rrr]^{\eta_{U,V,W,X}\otimes_1id_{Y\otimes_2Z}}
\ar[dd]^{id_{U\otimes_2V}\otimes_1\eta_{W,X,Y,Z}}&&&
\bigl((U\otimes_1W)\otimes_2(V\otimes_1X)\bigr)\otimes_1(Y\otimes_2Z)\ar[dd]^{\eta_{U\otimes_1W,V\otimes_1X,Y,Z}}\\
\\
(U\otimes_2V)\otimes_1\bigl((W\otimes_1Y)\otimes_2(X\otimes_1Z)\bigr)\ar[rrr]^{\eta_{U,V,W\otimes_1Y,X\otimes_1Z}}&&&
(U\otimes_1W\otimes_1Y)\otimes_2(V\otimes_1X\otimes_1Z)}
\end{equation}
The external associativity condition gives the commutative diagram:
\begin{equation}\label{extexpl}
\xymatrix{(U\otimes_2V\otimes_2W)\otimes_1(X\otimes_2Y\otimes_2Z)\ar[rrr]^{\eta_{U\otimes_2V,W,X\otimes_2Y,Z}}
\ar[dd]^{\eta_{U,V\otimes_2W,X,Y\otimes_2Z}}&&&\bigl((U\otimes_2V)\otimes_1(X\otimes_2Y)\bigr)\otimes_2(W\otimes_1Z)
\ar[dd]^{\eta_{U,V,X,Y}\otimes_2id_{W\otimes_1Z}}\\
\\
(U\otimes_1X)\otimes_2\bigl((V\otimes_2W)\otimes_1(Y\otimes_2Z)\bigr)\ar[rrr]^{id_{U\otimes_1X}\otimes_2\eta_{V,W,Y,Z}}&&&
(U\otimes_1X)\otimes_2(V\otimes_1Y)\otimes_2(W\otimes_1Z) }
\end{equation}

Finally, [BFSV] gives
\begin{defn}{\rm
Denote by ${\mathscr{MC}at}_n$ the category of (small) $n$-fold
monoidal categories. Then an $(n+1)$-fold monoidal category is a
monoid in ${\mathscr{MC}at}_n$.}
\end{defn}

This gives the following compatibility axiom: for $1\le i<j<k\le n$ the following diagram is commutative:

\begin{equation}\label{compexpl}
{\scriptsize\xymatrix{
&\bigl(((A_1\otimes_kA_2)\otimes_j(B_1\otimes_kB_2)\bigr)\otimes_i\bigl((C_1\otimes_kC_2)\otimes_j(D_1\otimes_kD_2)\bigr)
\ar[ldd]^{\eta^{jk}\otimes_i\eta^{jk}}\ar[dd]_{\eta^{ij}}\\
\\
\bigl((A_1\otimes_jB_1)\otimes_k(A_2\otimes_jB_2)\bigr)\otimes_i\bigl((C_1\otimes_jD_1)\otimes_k(C_2\otimes_j D_2)\bigr)\ar[dd]^{\eta^{ik}}&
\bigl((A_1\otimes_k A_2)\otimes_i(C_1\otimes_kC_2)\bigr)\otimes_j\bigl((B_1\otimes_kB_2)\otimes_i(D_1\otimes_kD_2)\bigr)\ar[dd]_{\eta^{ik}\otimes_j\eta^{ik}}\\
\\
\bigl((A_1\otimes_jB_1)\otimes_i(C_1\otimes_jD_1)\bigr)\otimes_k\bigl((A_2\otimes_jB_2)\otimes_i(C_2\otimes_jD_2)\bigr)\ar[ddr]^{\eta^{ij}\otimes_k\eta^{ij}}&
\bigl((A_1\otimes_iC_1)\otimes_k(A_2\otimes_iC_2)\bigr)\otimes_j\bigl((B_1\otimes_iD_1)\otimes_k(B_2\otimes_iD_2)\bigr)\ar[dd]_{\eta^{jk}}\\
\\
&\bigl((A_1\otimes_iC_1)\otimes_j(B_1\otimes_iD_1)\bigr)\otimes_k\bigl((A_2\otimes_iC_2)\otimes_j(B_2\otimes_iD_2)\bigr)
}}
\end{equation}

\subsection{\sc Examples}
\begin{example}\label{ex1}{\rm
Let $A$ be an associative algebra. The category $\Bimod(A)$ of $A$-bimodules is monoidal (1-fold monoidal), with $\otimes_A$ as the monoidal product.
The tautological bimodule $A$ is the unit in this monoidal category.
}
\end{example}

\begin{example}\label{ex2}{\rm
Let $B$ be an associative bialgebra over a field $k$ (see Section 2.1 below for definition). The category of left $B$-modules is monoidal, with the following monoidal structure. Let $M,N$ be two left $B$-modules, then $M\otimes_k N$ is naturally a $B\otimes_k B$-module. Now the coproduct $\Delta\colon B\to B\otimes_k B$ is a homomorphism of algebras; this makes $M\otimes_k N$ a $B$-module. Let $\varepsilon \colon B\to k$ be the counit in $B$, it endows $k$ with a $B$-module structure. This left $B$-module $k$ is the unit in this monoidal category.
}
\end{example}

\begin{example}\label{examplejs}{\rm
According to a result of Joyal and Street [JS], there are just few examples for $n\ge 2$ when the map $\eta^{ij}_{A,B,C,D}$ are {\it isomorphisms}
for any $A,B,C,D$ and any $1\le i<j\le n$ and when {\it there is a common unit object for all $n$ monoidal structures}. For $n=2$ any such category is equivalent as a 2-fold monoidal category to a category
with $A\otimes_1 B=A\otimes_2 B$ with a {\it braiding} $c_{A,B}\colon A\otimes B\to B\otimes A$ defining a structure of a {\it braided category} (see, e.g., [ES]) on
$\mathcal{C}$. Then we can construct a map $\eta_{A,B,C,D}\colon
(A\otimes B)\otimes(C\otimes D)\to (A\otimes C)\otimes (B\otimes D)$
just as $\eta_{A,B,C,D}=id\otimes c_{23}\otimes id$. This
construction gives a 2-fold monoidal category. For $n>2$ and $\eta_{A,B,C,D}$ isomorphisms one necessarily
has $A\otimes_i B=A\otimes_j B$ for any $i,j$ and {\it all $\otimes_i$ are symmetric}.
See the Fiedorowicz's Obervolfach talk [F] (page 4 and thereafter) for a clear overview of this result.}
\end{example}

\begin{example}\label{extetra}
{\rm
Let $A$ be an associative bialgebra. We define a {\it tetramodule}
over it as a $k$-vector space $M$ such that there is a bialgebra
structure on $A\oplus\epsilon M$, where $\epsilon^2=0$ and the
restriction of the bialgebra structure to $A$ is the initial one
(see Section 2 for details). If we rephrase this definition replacing
``bialgebra'' by ``associative algebra'', we recover the concept of a
bimodule over an associative algebra; thus, this Example is a generalization of Example \ref{ex1}. We
construct a 2-fold monoidal structure on the abelian
category $\Tetra(A)$ of tetramodules over $A$ in Section ???.}
\end{example}

\begin{example}\label{monalg}
{\rm
Examples \ref{ex1} and \ref{extetra} can be generalized as follows. Recall
from Example \ref{ex2} that the category of left modules over an associative bialgebra
is a monoidal category, with the monoidal structure equal to the
tensor product over $k$ on the level of the underlying vector spaces. Define an {\it
$n$-fold monoidal bialgebra} as an associative algebra with $n$
coassociative coproducts $\Delta_1,\dots,\Delta_n\colon A\to
A\otimes_k A$ such that the corresponding $n$ monoidal structures on
the category of left $A$-modules form an $n$-fold monoidal category.
Thus, 0-monoidal bialgebra is just an associative algebra, and
1-monoidal bialgebra is a bialgebra. One can define the category of
tetramodules over an $n$-monoidal bialgebra analogously to the
previous Example. We claim that this category is an $(n+1)$-fold
monoidal category; a proof is straightforward. We think this definition of $n$-monoidal bialgebra is a conceptually right
$n$-categorical generalization of the concept of bialgebra, for higher $n$.}
\end{example}

\subsection{\sc The operad of categories governing the $n$-fold monoidal
categories}
Fix $n\ge 1$. For any $d\ge 0$ denote by $\mathcal{M}_n(d)$ the full subcategory of the free $n$-fold monoidal category generated by objects $x_1,\dots, x_d$ consisting of objects which are monomials in $x_i$, where each $x_i$ occurs exactly ones.
For example, such monomials for $d=3$ and $n=2$ could be $(x_3\otimes_1x_1)\otimes_2 x_2$, or $(x_2\otimes_2 x_3)\otimes_1 x_1$. For fixed $n$ and $d$ the category $\mathcal{M}_n(d)$ has a finite number of objects. The morphisms in $\mathcal{M}_n(d)$ are exactly those which can be obtained as compositions of the associativities for a fixed $\otimes_i$, and $\eta_{ijkl}$, with exactly the same commutative diagrams as in $n$-fold monoidal category.

When $n$ is fixed and $d$ is varied, the categories $\mathcal{M}_n(d)$ form an operad of categories.
The following lemma follows from the definitions.

\begin{lemma}
A category is $n$-fold monoidal if and only if there is an action of the operad $\{\mathcal{M}_n(d)\}_{d\ge 0}$ of categories on it.
\end{lemma}
\endproof

The following very deep theorem is proven in [BFSV]:
\begin{theorem}
The classifying space of the operad of categories $\{\mathcal{M}_n(d)\}$ is an operad of topological space which is homotopically equivalent (as operad) to the $n$-dimensional little discs operad.
\end{theorem}

\section{\sc The category of tetramodules and a 2-fold monoidal structure on it}
\subsection{\sc }
Recall that an {\it associative bialgebra} is a vector space $A$ over a field $k$
equipped with two operations, the product $*:A^{\otimes 2}\to A$ and
the coproduct $\Delta\colon A\to A^{\otimes 2}$, which obey the
axioms 1.-4. below:
\begin{itemize}
\item[1.] Associativity: $a*(b*c)=(a*b)*c$;
\item[2.] Coassociativity: $(\Delta\otimes
id)\Delta(a)=(id\otimes\Delta)\Delta(a)$;
\item[3.] Compatibility: $\Delta(a*b)=\Delta(a)*\Delta(b)$.
\end{itemize}
We use the classical notation
$$
\Delta(a)=\Delta^1(a)\otimes\Delta^2(a)
$$
which is just a simplified form of the equation
$$
\Delta(a)=\sum_i\Delta^1_i(a)\otimes\Delta^2_i(a)
$$
We always assume that our bialgebras have a unit and a counit.
A unit is a map $i\colon k\to A$ and the counit is a map $\varepsilon\colon A\to k$. We always assume
\begin{itemize}
\item[4.] $i(k_1\cdot k_2)=i(k_1)*i(k_2)$, $\varepsilon(a*b)=\varepsilon(a)\cdot \varepsilon(b)$.
\end{itemize}
We also denote the product $*$ by $m$.

A {\it Hopf algebra} is a bialgebra with antipode. An antipode is a $k$-linear map $S\colon A\to A$ which obeys
\begin{itemize}
\item[5.] $m(1\otimes S)\Delta(a)=m(S\otimes 1)\Delta(a)=i(\varepsilon(a))$
\end{itemize}

As we already mentioned above, an ``operadic'' definition of a
bimodule over an associative algebra $A$ reads: it is a $k$-vector space $M$
such that $A\oplus\epsilon M$ is again an associative algebra, where
$\epsilon^2=0$, and the restriction of the algebra structure to $A$
coincides with the initial one. (The latter condition can be formulated equivalently that the associative algebra $A\oplus \epsilon M$ is defined {\it over} $A$). A vector space $M$ obeying this definition is the same that an
$A$-bimodule. We give an analogous definition in the case when $A$
is an associative bialgebra.

\begin{defn}{\rm
Let $A$ be an associative bialgebra. A {\it Bernstein-Khovanova
tetramodule} [BKh], [Kh] $M$ over $A$ is a vector space such that $A\oplus
\epsilon M$ is an associative bialgebra,  when $\epsilon^2=0$ and the
restriction of the bialgebra structure to $A$ is the initial one.
The category of tetramodules over a bialgebra $A$ is denoted
$\Tetra(A)$.}
\end{defn}

More precisely, one has maps $m_\ell\colon A\otimes M\to M$,
$m_r\colon M\otimes A\to M$ (which make $M$ an $A$-bimodule), and
maps $\Delta_\ell\colon M\to A\otimes M$ and $\Delta_r\colon M\to
M\otimes A$ (which make $M$ an $A$-bicomodule), with some
compatibility between these 4 maps. The compatibility written down
explicitly is the following 4 equations:
\begin{equation}\label{eq4.10}
\Delta_\ell(a*m)=(\Delta^1(a)*\Delta_\ell^1(m))\otimes
(\Delta^2(a)*\Delta^2_\ell(m))\subset A\otimes_k M
\end{equation}

\begin{equation}\label{eq4.11}
\Delta_\ell(m*a)=(\Delta_\ell^1(m)*\Delta^1(a))\otimes
(\Delta_\ell^2(m)*\Delta^2(a))\subset A\otimes_k M
\end{equation}

\begin{equation}\label{eq4.12}
\Delta_r(a*m)=(\Delta^1(a)*\Delta_r^1(m))\otimes
(\Delta^2(a)*\Delta^2_r(m))\subset  M\otimes_k A
\end{equation}

\begin{equation}\label{eq4.13}
\Delta_r(m*a)=(\Delta_r^1(m)*\Delta^1(a))\otimes
(\Delta^2_r(m)*\Delta^2(a))\subset M\otimes_k A
\end{equation}
Here we use the natural notation like
$\Delta_\ell(m)=\Delta^1_\ell(m)\otimes\Delta^2_\ell(m)$ with
$\Delta_\ell^1(m)\in A$, $\Delta^2_\ell(m)\in M$, etc. As well, we
use the sign $*$ for the both product in $A$ and the module products
$m_\ell$ and $m_r$.

The tetramodules over a bialgebra $A$ form an abelian category.
The main (non-trivial) example of a tetramodule over $A$ is $A$ itself; it is called
{\it the tautological tetramodule}.

When $A$ is finite-dimensional over $k$, a tetramodule is the same
that a left module over some associative algebra $H(A)$. This
algebra $H(A)$ is, as a vector space, the tensor product
$H(A)=A\otimes_kA\otimes_kA^*\otimes_k A^*$, and the product in it is defined to fulfill
the equations \eqref{eq4.10}-\eqref{eq4.13}
above. (An analogous associative algebra constructed from a bialgebra $A$, whose underlying vector space is $A\otimes_k A^*$, and the product is defined to fulfill only the single relation \eqref{eq4.10} among the 4 relations above, is known as {\it the Heisenberg double} of $A$.).

In particular, if $A$ is finite-dimensional over $k$, the abelian category
$\Tetra(A)$ has enough projective and enough injective objects.
For general $A$, R.Taillefer proved [Tai2] that the category
$\Tetra(A)$ has enough injectives.

Denote by $H_\GS^\mb(A,A)$ the Gerstenhaber-Schack cohomology [GS] of a bialgebra $A$ (we recall the definition in Section 3.1 below).
The Gerstenhaber-Schack cohomology $H^\mb_\GS(A,A)$ is known to
control the infinitesimal deformations of the bialgebra $A$.

The interplay between the category of tetramodules and the deformation theory of associative bialgebras is given in the following result:

\begin{theorem}[{Taillefer, [Tai1,2]}]\label{tt}
For any bialgebra $A$ one has:
\begin{equation}\label{eq4.14}
H^\mb_\GS(A,A)=\Ext^\mb_{\Tetra(A)}(A,A)
\end{equation}
\end{theorem}
\endproof

We recall the original Taillefer's proof of Theorem \ref{tt} in Section 3 below.

\begin{example}{\rm
Consider the case when $A=S(V)$ is a free (co)commutative bialgebra,
for $V$ a vector space over $k$. Suppose for simplicity that $V$ is finite-dimensional. We prove in Section 4 that
\begin{equation}\label{eq4.15}
H^k_\GS(A,A)=\oplus_{i+j=k}\Lambda^iV\otimes_k\Lambda^jV^*
\end{equation}
The total graded space $F=\oplus_{i,j\ge 0}\Lambda^iV\otimes_k\Lambda^jV^*[-i-j]$ can be thought as the space of functions on the space
$W=V[1]\oplus V^*[1]$. The space $W$ is a Poisson space, with the Poisson bracket of degree -2 given by the contraction of $V[1]$ with $V^*[1]$.
Altogether, $F$ admits a product and an even Poisson-Lie bracket of degree -2, which are compatible by the Leibniz rule. It can be rephrased saying that $F$ is a 3-algebra.

In fact, the Deligne conjecture for $2$-fold monoidal categories (mentioned in Introduction), and our Theorem \ref{mainresult} below, would imply together that the Gerstenhaber-Schack cohomology of any Hopf algebra is naturally a 3-algebra structure. Presumably, it is not true for general bialgebra which is not a Hopf algebra.
}
\end{example}

\subsection{\sc The structure of a 2-fold monoidal category on
$\Tetra(A)$}

\subsubsection{\sc Two ``external'' tensor products}
Before defining two ``real'' monoidal structures on the category of tetramodules, we start with two preliminary monoidal products. Our preliminary products are related to the real ones, as $M\otimes_k N$ is related to $M\otimes_A N$, for an associative algebra $A$ over $k$, and two its bimodules.

Let $A$ be an associative bialgebra, and let $M_1,M_2$ be two tetramodules over it.
We define two their
``external'' tensor products $M_1\boxtimes_1M_2$ and
$M_1\boxtimes_2M_2$ (which are $A$-tetramodules again). In the both
cases the underlying vector space is $M_1\otimes_kM_2$.

{The case of $M_1\boxtimes_1M_2$}:
\begin{itemize}
\item[1.] $m_\ell(a\otimes m_1\boxtimes m_2)=(am_1)\boxtimes m_2$,
\item[2.] $m_r(m_1\boxtimes m_2\otimes a)=m_1\boxtimes (m_2a)$,
\item[3.] $\Delta_\ell(m_1\boxtimes
m_2)=(\Delta_\ell^1(m_1)*\Delta_\ell^1(m_2))\otimes
(\Delta^2_\ell(m_1)\boxtimes\Delta^2_\ell(m_2))$,
\item[4.] $\Delta_r(m_1\boxtimes m_2)=(\Delta^1_r(m_1)\boxtimes
\Delta^1_r(m_2))\otimes (\Delta^2_r(m_1)*\Delta^2_r(m_2))$.
\end{itemize}

{The case of $M_1\boxtimes_2 M_2$}:
\begin{itemize}
\item[1.] $m_\ell(a\otimes m_1\boxtimes
m_2)=(\Delta^1(a)m_1)\boxtimes (\Delta^2(a)m_2)$,
\item[2.] $m_r(m_1\boxtimes m_2\otimes
a)=(m_1\Delta^1(a))\boxtimes(m_2\Delta^2(a))$,
\item[3.] $\Delta_\ell(m_1\boxtimes
m_2)=\Delta^1_\ell(m_1)\otimes(\Delta_\ell^2(m_1)\boxtimes m_2)$,
\item[4.] $\Delta_r(m_1\boxtimes
m_2)=(m_1\boxtimes\Delta_r^1(m_2))\otimes\Delta_r^2(m_2)$.
\end{itemize}

The main lack of the two external products is that the tautological tetramodule $A$ is not a unit object for them.

In the both definitions we use only ``the half'' of the tetramodule
structures on $M_1,M_2$. In particular, in the first definition we
do not use the right multiplication $m_r$ for $M_1$ and the left
multiplication $m_\ell$ for $M_2$. Similarly, in the second
definition we do not use $\Delta_r$ for $M_1$ and $\Delta_\ell$ for
$M_2$.

This provides us some extra possibilities to modify our two products, which we use below to define
the real ``internal'' tensor products $M_1\otimes_1 M_2$ and
$M_1\otimes_2 M_2$. For the two internal products, the tautological tetramodule $A$
is a unit object.

\subsubsection{\sc Two ``internal'' tensor products}
\begin{defn}{\rm
Let $M_1,M_2$ be two tetramodules over a bialgebra $A$. Their first
tensor product $M_1\otimes_1 M_1$ is defined as the
quotient-tetramodule
\begin{equation}\label{mon1}
M_1\otimes_1M_2=M_1\boxtimes_1M_2/((m_1a)\boxtimes_1m_2-m_1\boxtimes_1(am_2))
\end{equation}
One easily checks that this definition is correct. Analogously, the
second tensor product $M_1\otimes_2M_2$ is defined as a
sub-tetramodule
\begin{equation}\label{mon2}
M_1\otimes_2M_2=\left\{\sum_im_{1i}\boxtimes_2 m_{2i}\subset
M_1\boxtimes_2 M_2|\sum_i \Delta_r(m_{1i})\otimes_k
m_{2i}=\sum_im_{1i}\otimes_k\Delta_\ell(m_{2i})\right\}
\end{equation}
Again, one easily checks that this definition is correct.}
\end{defn}
\begin{lemma}\label{units}
Let $A$ be an associative bialgebra. Then the
tautological tetramodule $A$ is the unit object for both
monoidal structures.\end{lemma} \proof{} Let $M$ be a tetramodule.
One can check that the  maps
\begin{equation}\label{a1}
\begin{aligned}
\ & m_\ell\colon A\otimes_1M\to M\\
&m_r\colon M\otimes_1A\to M
\end{aligned}
\end{equation}
and
\begin{equation}\label{a2}
\begin{aligned}
\ & \Delta_\ell\colon M\to A\otimes_2M\\
&\Delta_r\colon M\to M\otimes_2A
\end{aligned}
\end{equation}
are morphisms of tetramodules.
As $A$ has a unit, the first two maps are bijective; as $A$ has a counit, the last two maps are bijective too.
\endproof

\subsubsection{\sc The 2-fold monoidal structure}
We construct the map $\eta_{M,N,P,Q}\colon (M\otimes_2
N)\otimes_1(P\otimes_2 Q)\to (M\otimes_1 P)\otimes_2(N\otimes_1 Q)$
in several steps.

\smallskip
\smallskip

{\it The first step} is to check that the map $\phi_0\colon
(M\boxtimes_2 N)\boxtimes_1(P\boxtimes_2 Q)\to (M\boxtimes_1
P)\boxtimes_2 (N\boxtimes_1 Q)$, $\phi_0(m\otimes_k n\otimes_k
p\otimes_k q)=m\otimes_k p\otimes_k n\otimes_k q$, is a map of
tetramodules. We have:
\begin{equation}\label{eqformula1}
a*\bigl((m\boxtimes_2 n)\boxtimes_1 (p\boxtimes_2
q)\bigr)=
\bigl(a*(m\boxtimes_2 n)\bigr)\boxtimes_1 (p\boxtimes_2
q)=
(\Delta^1(a)*m)\otimes(\Delta^2(a)*n)\otimes p\otimes q
\end{equation}
and
\begin{equation}\label{eqformula2}
a*\bigl((m\boxtimes_1 p)\boxtimes_2(n\boxtimes_1 q)\bigr)=
\bigl(\Delta^1(a)*(m\boxtimes_1
p)\bigr)\boxtimes_2\bigl(\Delta^2(a)*(n\boxtimes_1 q)\bigr)=
(\Delta^1(a)*m)\otimes p\otimes (\Delta^2(a)*n)\otimes q
\end{equation}
We see that
\begin{equation}
\phi_0(\text{r.h.s. of \eqref{eqformula1}})=\text{(r.h.s. of
\eqref{eqformula2})}
\end{equation}
That is, $\phi_0$ is a map of left modules; analogously it is a map
of right modules.

Now prove that $\phi_0$ is a map of left comodules. We have:
\begin{equation}\label{eqformula3}
\begin{aligned}
\ &\Delta_\ell\bigl((m\boxtimes_2 n)\boxtimes_1 (p\boxtimes_2
q)\bigr)=
\bigl(\Delta_\ell^1(m\boxtimes_2 n)*\Delta_\ell^1(p\boxtimes_2
q)\bigr)\otimes_k\bigl(\Delta_\ell^2(m\boxtimes_2
n)*\Delta_\ell^2(p\boxtimes_2 q)\bigr)=\\
&\bigl(\Delta_\ell^1(m)*\Delta_\ell^1(p)\bigr)\otimes_k\bigl(\Delta_\ell^2(m)\otimes_k
n\otimes_k\Delta_\ell^2(p)\otimes_k q\bigr)
\end{aligned}
\end{equation}
and
\begin{equation}\label{eqformula4}
\begin{aligned}
\ &\Delta_\ell\bigl((m\boxtimes_1 p)\boxtimes_2(n\boxtimes_1
q)\bigr)=
\Delta_\ell^1(m\boxtimes_1p)\otimes_k\bigl(\Delta_\ell^2(m\boxtimes_1
p)\boxtimes_2(n\boxtimes_1 q)\bigr)=\\
&\bigl(\Delta_\ell^1(m)*\Delta_\ell^1(p)\bigr)\otimes_k\bigl(\Delta_\ell^2(m)\otimes_k
\Delta_\ell^2(p)\otimes_k n\otimes_k q\bigr)
\end{aligned}
\end{equation}
We see that
\begin{equation}
\Delta_\ell\circ \phi_0=id\otimes_k(\phi_0\circ\Delta_\ell)
\end{equation}
that is, $\phi_0$ is a map of left comodules. It is proven
analogously that $\phi_0$ is a map of right comodules.

\vspace{2mm}

{\it At the second step} we consider the natural projections of
tetramodules $p_{M,P}\colon M\boxtimes_1 P\to M\otimes_1 P$ and
$p_{N,Q}\colon N\boxtimes_1Q\to N\otimes_1 Q$. We consider the
composition
\begin{equation}
\phi_1=(p_{M,P}\boxtimes_2 p_{N,Q})\circ\phi_0\colon (M\boxtimes_2
N)\boxtimes_1(P\boxtimes_2 Q)\to (M\otimes_1
P)\boxtimes_2(N\otimes_1 Q)
\end{equation}
We want to check that the map $\phi_1$ defines naturally a map
\begin{equation}
\phi_2=\overline{\phi}_1\colon (M\boxtimes_2
N)\otimes_1(P\boxtimes_2 Q)\to (M\otimes_1 P)\boxtimes_2(N\otimes_1
Q)
\end{equation}
that is, the elements of the form
\begin{equation}\label{eqformula10}
\bigl((m\boxtimes_2 n)*a\bigr)\otimes_k (p\boxtimes_2
q)-(m\boxtimes_2 n)\otimes_k \bigl(a*(p\boxtimes_2 q)\bigr)
\end{equation}
are mapped to 0 by $\phi_1$.

Indeed,
\begin{equation}
\eqref{eqformula10}=(m*\Delta^1a)\otimes_k(n*\Delta^2a)\otimes_k(p\otimes_k
q)-(m\otimes_k n)\otimes_k (\Delta^1a*p)\otimes_k (\Delta^2a*q)
\end{equation}
which, after the permutation $\phi_0$ of the second and the third
factors, is mapped to 0 in $(M\otimes_1 P)\boxtimes_2(N\otimes_1
Q)$. Therefore, the map $\phi_2$ is well-defined.

\vspace{2mm}

{\it At the third step} we restrict the map $\phi_2$ to $(M\otimes_2
N)\otimes_1(P\otimes_2 Q)\subset (M\boxtimes_2 N)\otimes_1
(P\boxtimes_2 Q)$, and we need to check that the image of this
restricted map belongs to $(M\otimes_1 P)\otimes_2(N\otimes_1
Q)\subset (M\otimes_1 P)\boxtimes_2(N\otimes_1 Q)$.

Suppose $m\otimes_k n\in M\otimes_2N$ and $p\otimes_k q\in
P\otimes_2 Q$ (we assume the summation over several such monomials,
but for simplicity we skip this summation). Then
\begin{equation}\label{eqformula11}
\Delta_r^1(m)\otimes_k\Delta_r^2(m)\otimes_k
n=m\otimes_k\Delta_\ell^1(n)\otimes_k\Delta_\ell^2(n)
\end{equation}
with the middle factors in $A$, and analogously
\begin{equation}\label{eqformula12}
\Delta_r^1(p)\otimes_k\Delta_r^2(p)\otimes_kq=p\otimes_k\Delta_\ell^1(q)\otimes_k\Delta_\ell^2(q)
\end{equation}
again, with the middle factors in $A$.

One needs to prove that \eqref{eqformula11} and \eqref{eqformula12}
together imply that
\begin{equation}\label{eqformula13}
(m\boxtimes_1 p)\boxtimes_2(n\boxtimes_1 q)\in (M\otimes_1
P)\otimes_2 (N\otimes_1 Q)\subset(M\otimes_1 P)\boxtimes_2
(N\otimes_1 Q)
\end{equation}
that is,
\begin{equation}\label{eqformula14}
\Delta_r^1(m)\otimes_k
\Delta_r^2(p)\otimes_k(\Delta_r^2(m)*\Delta_r^2(p))\otimes_k
n\otimes_kq=m\otimes_k p\otimes_k
(\Delta_\ell^1(n)*\Delta_\ell^1(q))\otimes_k
\Delta_\ell^2(n)\otimes_k\Delta_\ell^2(q)
\end{equation}

To get \eqref{eqformula14} from \eqref{eqformula12} and
\eqref{eqformula13} we permute \eqref{eqformula12} such that the
factors in $A$ are the most right, permute \eqref{eqformula12} such
that the factors in $A$ are the most left, then take the equation
$(\text{l.h.s. of }\eqref{eqformula12})\otimes_k(\text{l.h.s. of
}\eqref{eqformula13})=(\text{r.h.s. of
}\eqref{eqformula12})\otimes_k(\text{r.h.s. of
}\eqref{eqformula13})$ (after they are permuted). Then for the two
middle factor (in $A$) we apply the product $*\colon A\otimes_k A\to
A$, and then permute again.

\vspace{3mm}

{\it The map $\eta_{M,N,P,Q}$ is constructed.} \vspace{2mm}

\begin{theorem}\label{mainresult}
The maps $\eta_{M,N,P,Q}$, constructed above, and the two
tensor products $\otimes_1$ and $\otimes_2$, define a 2-fold
monoidal structure on the category $\Tetra(A)$ of tetramodules over
a bialgebra $A$. The tautological tetramodule $A$ is a unit object of this 2-fold monoidal category.
\end{theorem}
\proof{} First of all, the two tensor products $\boxtimes_1$ and
$\boxtimes_2$ with $\wtilde{\eta}_{M,N,P,Q}=\phi_0$, clearly define a
2-fold monoidal structure (without a unit object) on the category $\Tetra(A)$. In
particular, the diagrams \eqref{intexpl} and \eqref{extexpl} are
commutative for $\wtilde{\eta}_{M,N,P,Q}$ (because $\phi_0$ is just
the permutation which switches the second and the third factors).
Now the same diagrams for the actual structure $\eta_{M,N,P,Q}$ are
obtained from the above ``external structure'' diagrams just by passing to subquotients.
Therefore, they are commutative as well.
\endproof

\begin{remark}{\rm
The ``external'' 2-fold monoidal structure does not obey properly the definition of $n$-fold monoidal categories given in Section 1, as it does not have a unit object (one easily sees that $k$ is not a unit). Therefore, the Joyal-Street result from Example \ref{examplejs} is not applied to it.
}
\end{remark}

\section{\sc Tetramodules and the Gerstenhaber-Schack cohomology: an overview of Taillefer's results}
In this Section we present some results of R.Taillefer [Tai1,2] on the category of tetramodules. 
Our presentation follows the same ideas as the original proofs of Taillefer,
but the details are slightly different, and, we hope, more transparent.
In particular, we prove Theorem \ref{tt} for any bialgebra, whereas Taillefer assumes
it to be a Hopf algebra.

As an application of the general theory, we present a computation of the Gerstenhaber-Schack cohomology for $A=S(V)$, the free commutative cocommutative bialgebra generated by a vector space $V$.

\subsection{\sc The Gerstenhaber-Schack complex}
Let $A$ be a (co)associative bialgebra.  Note that the
bar-differential in $\BBar^{\boxtimes_1}(A)$ is given by maps of
tetramodules; analogously, the cobar-differential in
$\Cobar^{\boxtimes_2}(A)$ is given by maps of tetramodules.

Let us recall, that originally the Gerstenhaber-Schack complex was
defined in [GS] as

\begin{equation}\label{gs}
C_{\GS}^\mb(A)=\Hom_{\Tetra(A)}(\BBar^{\boxtimes_1}_-(A),\Cobar^{\boxtimes_2}_+(A))
\end{equation}

Here $\BBar_-(B)$ and $\Cobar_+(C)$ are truncated complexes, which
end (start) with $B{\boxtimes_1}B$ ($C{\boxtimes_2}C$)
correspondingly.

For convenience of the reader let us write down here the
Gerstenhaber-Schack differential in $C_{\GS}^\mb(A)$ explicitly:

First of all, as a graded vector space,
\begin{equation}\label{eqeqeq_1}
C_{\GS}^\mb(A)=\oplus_{m,n\ge 0}\Hom_k(A^{\otimes m}, A^{\otimes
n})[-m-n]
\end{equation}
Now let $\Psi\colon A^{\otimes m}\to A^{\otimes n}\in
C_{GS}^{m+n}(A)$. We are going to define the Gerstenhaber-Schack
differential $d_\GS (\Psi)\in \Hom (A^{\otimes (m+1)},A^{\otimes
n})\oplus\Hom (A^{\otimes m},A^{\otimes (n+1)})$. Denote the
projection of $d_\GS$ to the first summand by $(d_\GS)_1$, and the
projection to the second summand by $(d_\GS)_2$. The formulas for
$(d_\GS)_1$ and $(d_\GS)_2$ are:
\begin{equation}\label{eq0.6fin}
\begin{aligned}
\ &(d_\GS)_1(\Psi)(a_0\otimes\dots\otimes a_m)=\\
&\Delta^{n-1}(a_0)*\Psi(a_1\otimes\dots\otimes a_m)+
\sum_{i=0}^{m-1}(-1)^{i+1}\Psi(a_0\otimes\dots\otimes (a_i*
a_{i+1})\otimes
\dots\otimes a_m)+\\
&(-1)^{m-1}\Psi(a_0\otimes\dots\otimes a_{m-1})*\Delta^{n-1}(a_m)
\end{aligned}
\end{equation}
and
\begin{equation}\label{eq0.7fin}
\begin{aligned}
\ &(d_\GS)_2(\Psi)(a_1\otimes\dots\otimes a_m)=\\
&(\Delta^{(1)}(a_1)*\Delta^{(1)}(a_2)*\dots
*\Delta^{(1)}(a_m))\otimes
\Psi(\Delta^{(2)}(a_1)\otimes\dots\otimes \Delta^{(2)}(a_m))+\\
&\sum_{i=1}^n(-1)^i\Delta_{i}\Psi(a_1\otimes\dots\otimes a_m)+\\
&(-1)^{n+1}\Psi(\Delta^{(1)}(a_1)\otimes\Delta^{(1)}(a_2)\otimes\dots
\otimes\Delta^{(1)}(a_m))\otimes(\Delta^{(2)}(a_1)*\Delta^{(2)}(a_2)*\dots
*\Delta^{(2)}(a_m))
\end{aligned}
\end{equation}

The goal of this Section is to recall the original proof, with some minor improvements in presentation, of the following result due to R.Taillefer [Tai1,2]:

\begin{theorem}[Taillefer]
For any bialgebra $A$ one has:
$$
\Ext^\mb_{\Tetra(A)}(A,A)=H^\mb\bigl(\Hom_{\Tetra(A)}(\BBar^{\boxtimes_1}_{-}(A),\Cobar^{\boxtimes_2}_{+}(A))\bigr)
$$
\end{theorem}

\subsection{\sc Two forgetful functors and their adjoints}
\subsubsection{\sc }
Let $A$ be a (co)associative bialgebra. Besides the category
$\Tetra(A)$, we can consider the categories $\Bimod(A)$ of
$A$-bimodules (when we consider $A$ as an algebra) and $\Bicomod(A)$
of $A$-bicomodules (when we consider $A$ as a coalgebra). Clearly
there are two exact forgetful functors $F_1\colon \Tetra(A)\to
\Bicomod(A)$ and $F_2\colon \Tetra(A)\to \Bimod(A)$. We have the
following
\begin{lemma}
Let $A$ be a bialgebra which has unit and counit. Then the functor
$F_1$ admits a left adjoint $L$ and the functor $F_2$ admits a right
adjoint $R$. The functors $L$ and $R$ are exact.
\end{lemma}
\begin{defn}{\rm
Let $A$ be a bialgebra, and let $N$ be an $A$-bicomodule, and let
$M$ be an $A$-bimodule.  The tetramodule $L(N)$ is called {\it an
induced tetramodule} (from $N$), and the tetramodule $R(M)$ is called
{\it a coinduced tetramodule} (from $M$). The induced and coinduced
tetramodules form full additive subcategories in the abelian
category $\Tetra(A)$. We denote them $\Tetra_\Ind(A)$ and
$\Tetra_\Coind(A)$, respectively.}
\end{defn}

{\it Proof of Lemma:} we set
\begin{equation}\label{adjoint1}
L(N)=A\boxtimes_1 N\boxtimes_1 A
\end{equation}
and
\begin{equation}\label{adjoint2}
R(M)=A\boxtimes_2 M\boxtimes_2 A
\end{equation}
(see Section 2.2.1 for the definitions of $\boxtimes_1$ and
$\boxtimes_2$). Rigorously, to write down formulas like this, $M$
and $N$ should be tetramodules. However, the definition of $M_1\boxtimes_1 M_2$ does not use the
right $A$-module structure in $M_1$ and the left $A$-module
structure in $M_2$. Analogously, in the definition of
$M_1\boxtimes_2 M_2$ does not use the right comodule structure in
$M_1$ and the left comodule structure in $M_2$. Therefore,
(\ref{adjoint1}) and (\ref{adjoint2}) make sense.

The adjunction properties of $L$ and $R$ are clear.

The exactness of $L$ and $R$ is clear from the constructions
(\ref{adjoint1}) and (\ref{adjoint2}).
\endproof

\subsubsection{\sc }
\begin{lemma}
Let $A$ be an associative bialgebra. Then any tetramodule
$M\in\Tetra(A)$ can be imbedded into a coinduced tetramodule, and
there is a surjection to $M$ from an induced tetramodule.
\end{lemma}
\proof{} Let $M$ be an $A$-tetramodule. Consider
$P(M)=A\boxtimes_1 M\boxtimes_1 A$, it is induced from the
bicomodule $F_1(M)$. The map $p\colon A\boxtimes_1 M\boxtimes_1 A\to
M$, $a\boxtimes_1m\boxtimes_1b\mapsto a\cdot m\cdot b$ is clearly a
map (and an epimorphism, because $A$ contains a unit) of
tetramodules. Analogously, the tetramodule $Q(M)=A\boxtimes_2
M\boxtimes_2 A$ is coinduced from the bimodule $F_2(M)$, and we have
a monomorphism $j\colon M\to A\boxtimes_2 M\boxtimes_2 A$, $m\mapsto
\Delta_\ell\circ\Delta_r(m)$.
\endproof

\begin{coroll}[{[Tai2]}]
For any bialgebra $A$ the category $\Tetra(A)$ has enough
injectives. \end{coroll}
\proof{} The functor $R$ is a right adjoint to
an exact functor $F_2$, and, therefore, maps injective objects to
injective (see [W], Prop. 2.3.10). Moreover, $R$ is left exact ([W],
Section 2.6). Therefore, it is sufficient to imbed $M$ as an
$A$-bimodule into an injective $A$-bimodule $I$ (it is a classical construction, see e.g. [W], Section
2.3). Then we apply the functor $R$ to this this imbedding of
$A$-bimodules. It gives an imbedding $j\colon M\to Q(M)$.
\endproof

The main fact about the induced and the coinduced tetramodules is
the following:

\begin{prop}
Let $A$ be an associative bialgebra. Then the functor
$X\mapsto \Hom_{\Tetra(A)}(X,Q)$ for fixed $Q\in\Tetra_{\Coind}(A)$
is an exact functor from $\Tetra_\Ind(A)^\opp$ to $\mathscr{V}ect$. As well,
the functor $Y\mapsto\Hom_{\Tetra(A)}(P,Y)$ for fixed
$P\in\Tetra_\Ind(A)$ is an exact functor from $\Tetra_\Coind(A)$ to
$\mathscr{V}ect$. \end{prop}
\proof{} We prove the first statement. Let
\begin{equation}\label{adjoint5}
0\rightarrow LN^\prime\rightarrow LN\rightarrow
LN^{\prime\prime}\rightarrow 0
\end{equation}
be an exact sequence of tetramodules,
$N,N^{\prime},N^{\prime\prime}\in \Bicomod(A)$. We need to prove that
the sequence
\begin{equation}\label{adjoint6}
0\rightarrow\Hom_{\Tetra(A)}(LN^{\prime\prime},RM)\rightarrow\Hom_{\Tetra(A)}(LN,RM)\rightarrow\Hom_{\Tetra(A)}(LN^\prime,RM)\rightarrow
0
\end{equation}
is exact for any $M\in \Bimod(A)$.

By the adjunction, the exactness of (\ref{adjoint6}) is equivalent
to the exactness of the sequence
\begin{equation}\label{adjoint7}
0\rightarrow\Hom_{\Bimod(A)}(F_2LN^{\prime\prime},M)\rightarrow\Hom_{\Bimod(A)}(F_2LN,M)\rightarrow\Hom_{\Bimod(A)}(F_2LN^\prime,M)\rightarrow
0
\end{equation}
The latter sequence is indeed exact as for any $N\in\Bicomod(A)$ the
$A$-bimodule $F_2LN$ is free and, therefore, projective.

The second statement is proven analogously.
\endproof

\subsection{\sc Some homological algebra}
We recall here some construction of homological algebra, based on the Grothendieck's interpretation of the derived functors as ``universal $\delta$-functors'' [Tohoku, 2.1-2.2]. This construction provides a useful way of computation of $\Ext$'s functors in abelian categories, using resolutions from ``semi-projective'' and ``semi-injective'' objects.
\subsubsection{\sc A $(\mathcal{P},\mathcal{Q})$-pair}
\begin{defn}\label{pq}{\rm
Let $\mathcal{A}$ be an abelian category, and let $\mathcal{P}$,
$\mathcal{Q}$ be two additive subcategories of $\mathcal{A}$. We say that $\mathcal{P}$ and $\mathcal{Q}$ form
a $({P},{Q})$-pair, if the following conditions are
fulfilled:
\begin{itemize}
\item[1.] the functor $\Hom(?,Q)$ is exact on $\mathcal{P}^\opp$ for
any $Q\in\mathcal{Q}$;
\item[2.] the functor $\Hom(P,?)$ is exact on $\mathcal{Q}$ for any
$P\in\mathcal{P}$;
\item[3.] for any object $M\in\mathcal{A}$, there is an epimorphism
$P\to M$ for $P\in\mathcal{P}$, and there is a monomorphism $M\to Q$
for $Q\in\mathcal{Q}$;
\item[4.] a stronger version of 3: for any short exact sequence
$0\to M_1\to M_2\to M_3\to 0$ in $\mathcal{A}$ the epimorphisms
$P_i\to M_i$, $p_i\colon P_i\in\mathcal{P}$, can be chosen such that
there is a map of complexes
\begin{equation}\label{delta}
\xymatrix{
0\ar[r]&P_1\ar[r]\ar[d]^{p_i}&P_2\ar[r]\ar[d]^{p_2}&P_3\ar[r]\ar[d]^{p_3}&0\\
0\ar[r]&M_1\ar[r]&M_2\ar[r]&M_3\ar[r]&0}
\end{equation}
where the upper line is an exact sequence. As well, also the dual condition
for $\mathcal{Q}$ for monomorphisms $q_i\colon M_i\to Q_i$,
$Q_i\in\mathcal{Q}$ is required.
\end{itemize}
}
\end{defn}
Note that the third condition guarantees that each object
$M\in\mathcal{A}$ has a $\mathbb{Z}_{\le 0}$-graded resolution in
$\mathcal{P}$ and a $\mathbb{Z}_{\ge 0}$-graded resolution in
$\mathcal{Q}$.
\begin{example}{\rm
If $\mathcal{A}$ has enough projectives and $\mathcal{P}$ is the
additive subcategory of projective objects,
$\mathcal{Q}=\mathcal{A}$ gives a $({P},{Q})$-pair.
Analogously, if $\mathcal{A}$ has enough injectives and
$\mathcal{Q}$ is the additive subcategory of injective objects,
$\mathcal{P}=\mathcal{A}$ gives a $({P},{Q})$-pair.}
\end{example}
\begin{prop}
Let $A$ be an associative bialgebra, and let
$\mathcal{A}=\Tetra(A)$. Then the pair
$(\Tetra_\Ind(A),\Tetra_\Coind(A))$ is a
$({P},{Q})$-pair.
\end{prop}
\proof{} The first two properties were proven in Proposition 3.6, and
the third property follows from Lemma 3.4. Moreover, the
construction in the Lemma gives immediately the fourth assertion in
Definition 3.7.
\endproof

\subsubsection{\sc The Key-Lemma}
The main fact about $({P},{Q})$-pairs is the
following lemma:
\begin{klemma}
Let $\mathcal{A}$ be an abelian category, having enough projective
or injective objects. Suppose we are given a
$({P},{Q})$-pair $(\mathcal{P},\mathcal{Q})$ in $\mathcal{A}$, and let $M,N\in\mathcal{A}$ be two
objects. Suppose $P^\mb\to M$ is a resolution of $M$ of objects in
$\mathcal{P}$, and $N\to Q^\mb$ is a resolution of $N$ by objects in
$\mathcal{Q}$. Then
\begin{equation}\label{kl}
\Ext^\mb_{\mathcal{A}}(M,N)=H^\mb(\Hom_{\mathcal{A}}(P^\mb,Q^\mb))
\end{equation}
\end{klemma}
\proof{} The proof consists from several steps. We give the proof
for the case when $\mathcal{A}$ has enough injectives, the case of enough projectives is
analogous.

The rough idea is: any derived functor in an abelian category (if it exists) enjoys the property of being of universal $\delta$-functor, and as such, is uniquely defined up to an isomorphism by its zero degree component. We prove that the functor $(M,N)\mapsto
H^\mb(\Hom_{\mathcal{A}}(P^\mb,Q^\mb))$ is a $\delta$-functor, and then prove that it is a universal $\delta$-functor.

Let $\mathcal{A}$ and $\mathcal{B}$ be two abelian categories, and let  $\{T_n\colon\mathcal{A}\to\mathcal{B}\}$, $n\ge 0$
be a collection of additive functors. One says that this collection is
{\it a cohomological $\delta$-functor} if for any exact sequence
\begin{equation}\label{ses}
0\to M\to N\to L\to 0
\end{equation}
in $\mathcal{A}$, one has a morphism $\delta\colon T_n(L)\to
T_{n+1}(M)$, $n\ge 0$ in $\mathcal{B}$, forming a long
exact sequence:
\begin{equation}\label{les}
\dots\rightarrow T_{n-1}(L)\xrightarrow{\delta}T_n(M)\rightarrow
T_n(N)\rightarrow T_n(L)\xrightarrow{\delta}
T_{n+1}(M)\rightarrow\dots
\end{equation}
$n\ge 1$, depending functorially on the short exact sequence
(\ref{ses}).

Consider a $\delta$-functor $\{T_n\}$. We say that this
$\delta$-functor is {\it universal} if for any other
$\delta$-functor $\{S_n\}$ with the natural transformation
$f_0\colon T_0\to S_0$ there is a unique morphism of
$\delta$-functors $\{f_n\colon T_n\to S_n\}$ extending $f_0$. It follows immediately that the universal $\delta$-functor with
$T_0=F$, if it exists, is unique. This point of view, independent on
existence of enough projective objects, was emphasized by Grothendieck in [Tohoku]
(see also [W], Chapter 2).

Now consider the functor $\Hom_{\mathcal{A}}(M,?)$ as a functor of the
second argument. If $\mathcal{A}$ has enough injectives, the functors
$T_n(M,N)=\Ext^n_{\mathcal{A}}(M,N)$ is a universal cohomological
$\delta$-functor (see [W], Theorem 2.4.7).

Now the proof of Key-Lemma goes in several steps:
\begin{itemize}
\item[Step 1.] $T_k\colon(M,N)\mapsto H^k(\Hom_{\mathcal{A}}(P^\mb,Q^\mb))$
with $P^\mb\in\mathcal{P}$, $Q^\mb\in\mathcal{Q}$ is well-defined,
that is does not depend on the choice of $P^\mb$ and $Q^\mb$;
\item[Step 2.] it is a cohomological $\delta$-functor with
$T_0=\Hom_{\mathcal{A}}(M,N)$;
\item[Step 3.] it is a universal cohomological $\delta$-functor.
\end{itemize}
Clearly the Key-Lemma follows from these three claims.

{\it Step 1:} it easily follows from the conditions 1. and 2. in the
definition of a $({P},{Q})$-pair.

{\it Step 2:} it follows easily from condition 4. in the definition
of a $(\mathcal{P},\mathcal{Q})$-pair.

{\it Step 3:} this is a bit more tricky. An additive functor
$F\colon \mathcal{A}\to\mathcal{B}$ is called {\it effaceable} if
for any object $N\in\mathcal{A}$ there is a monomorphism $j\colon
M\to I$ such that $F(j)=0$. It is proven in [Tohoku] that a
cohomological $\delta$-functor $\{T_n\}$ for which all $T_n$ for
$n\ge 1$ are effaceable, is universal. It remains to prove that our
functors $T_n(N)=H^n(\Hom_{\mathcal{A}}(P^\mb(M),Q^\mb(N)))$, $n\ge
1$, are effaceable. We can choose a monomorphism $j\colon N\to I$
with $I\in\mathcal{Q}$ by condition 3. in Definition \ref{pq}. Now the effaceability follows as
$T_n(I)=0$, $n\ge 1$, by 1. and 2. in Definition \ref{pq}.

Thus, it is proven that the functors
$H^n(\Hom_{\mathcal{A}}(P^\mb,Q^\mb))$ form a universal cohomological
$\delta$-functor which has the same the same 0-component as
$\Ext^n_{\mathcal{A}}(M,N)$. Therefore, the two $\delta$-functors are isomorphic.

\endproof

Theorem 3.1 is proven.
\endproof

\subsection{\sc Example: a computation of the Gerstenhaber-Schack
cohomology for $A=S(V)$} Here we compute, as an application of
previous results, the Gerstenhaber-Schack cohomology
for the (co)free commutative and cocommutative bialgebra $A=S(V)$.
For simplicity, we assume $V$ is a finite-dimensional, although the result (being slightly modified) remains true for general $V$.
\begin{prop}
Let $A=S(V)$ be (co)free commutative cocommutative bialgebra, $V$
finite-dimensional. Then the Gerstenhaber-Schack cohomology of $A$ is
$$H^k_\GS(A)=\bigoplus_{\substack{i+j=k\\ i,j\ge
0}}\Lambda^iV\otimes\Lambda^j(V^*)$$ \end{prop}

\proof{}
We compute $\Ext^\mb_{\Tetra(A)}(A,A)$ for $A=S(V)$, as $H^\mb(\Hom(P^\mb, Q^\mb))$, where $P^\mb$ is a resolution of $A$ by tetramodules which are free as bimodules, and $Q^\mb$ is a resolution of $A$ by tetramodules which are cofree as bicomodules. We are allowed to use such resolutions by Proposition 3.9 and Key-lemma 3.10.

We now use more ``economic'' resolutions $P^\mb$ and $Q^\mb$, than the bar- and cobar-resolutions from Theorem 3.1. Namely, we use suitable Koszul resolutions.

Recall these resolutions $P^\mb$ and $Q^\mb$. We have:

\begin{equation}
\begin{aligned}
\ &P^{-k}=S(V)\boxtimes_1 \Lambda^k(V)\boxtimes_1 S(V),\\
&\partial_k(f\otimes (e_0\wedge\dots\wedge e_k)\otimes g)=\\
&\sum_{i=0}^k(-1)^i\left((v_i\cdot f)\otimes (e_0\wedge\dots\wedge\hat{e_i}\wedge\dots\wedge e_k)\otimes g-f\otimes (e_0\wedge\dots\wedge\hat{e_i}\wedge\dots\wedge e_k)\otimes (v_i\cdot g)\right)
\end{aligned}
\end{equation}

\begin{equation}
\begin{aligned}
\ &Q^k=S(V)\boxtimes_2 \Lambda^k(V)\boxtimes_2 S(V),\\
&d_k(f\otimes (e_0\wedge\dots\wedge e_k)\otimes g)=\\
&\Delta_r^1(f) \otimes (\Delta_r^2(f)\wedge e_0\wedge\dots \wedge e_k)\otimes g-f\otimes(e_0\wedge\dots\wedge e_k\wedge \Delta_\ell^1(g))\otimes \Delta_\ell^2(g)
\end{aligned}
\end{equation}
where in the both equations the tetramodule structure on $\Lambda^k(V)$ is trivial for any $k$. In the second equation $\Delta_\ell(f)$ denotes
the projection of $\Delta(f)\in S(V)\otimes S(V)$ to $V\otimes S(V)$, and $\Delta_r(f)$ denotes the projection of $\Delta(f)$ to $S(V)\otimes V$.

One easily sees, that when $\Hom_{\Tetra(S(V))}(P^\mb, Q^\mb)$ is computed, the differential in this complex vanishes. This gives the result.

\endproof

\bigskip

\noindent {\sc Max-Planck Institut f\"{u}r Mathematik, Vivatsgasse 7, 53111 Bonn,\\
GERMANY}

\bigskip

\noindent {\em E-mail address\/}: {\tt borya$\_$port@yahoo.com}

\end{document}